\documentclass{siamltex}		
\usepackage{graphicx,multirow,amsfonts,subfigure,hyperref,color}
\usepackage[normalem]{ulem}
\usepackage[export]{adjustbox}
\usepackage{amssymb,amsmath}
      
\usepackage[utf8]{inputenc}

\usepackage{algorithm}
\usepackage{algpseudocode}
\def\ubarQ{{\underline{Q}}}
\def\ubarLambda{{\underline{\Lambda}}}
\def\half{\frac{1}{2}}

\newcommand{\eq}[1]{\begin{equation}\label{#1}}
\newcommand{\en}{\end{equation}}
\def\Off{\texttt{Off}} 
\def\Diag{\texttt{Diag}} 
\def\inv{^{-1}}%
\def\nref#1{(\ref{#1})}
\def\RR{\mathbb{R}}
\def\trace{\text{Trace}}

\def\up#1{^{(#1)}}%
\def\Span#1{\text{Span}{\{#1\}}}%

\graphicspath{{./FIGS/}}

\title{Joint Approximate Partial Diagonalization of Large Matrices}
\author{Abd-Krim~Seghouane\thanks{
Department of Electrical and Electronic Engineering, Melbourne School of 
Engineering,
The University of Melbourne, Melbourne, Australia;
E-mail: Abd-krim.seghouane@unimelb.edu.au}
\and Yousef Saad\thanks{
Department of Computer Science and Engineering,
The University of Minnesota at Twin Cities, Minnesota, USA. 
Email: saad@umn.edu.
Work of this author was supported by NSF  grant
CCF-1318597.}}

\date{\today} 

\begin{document} 

\maketitle 

\begin{abstract}
Given a  set of  $p$ symmetric (real)  matrices, the  Orthogonal Joint
Diagonalization (OJD) problem consists of finding an orthonormal basis
in which the representation of each  of these $p$ matrices is as close
as possible to a diagonal matrix.  We argue that when the matrices are
of large dimension, then the natural generalization of this problem is
to seek  an orthonormal  basis of  a certain subspace  that is  a near
eigenspace for all the  matrices in the set.  We refer  to this as the
problem of ``partial joint diagonalization of matrices.'' The approach
proposed first finds this approximate  common near eigenspace and then
proceeds to a  joint diagonalization of the restrictions  of the input
matrices in  this subspace.  A  few solution methods for  this problem
are  proposed  and illustrations  of  its  potential applications  are
provided.
\end{abstract}



\textbf{AMS:} 15A69, 15A18


\section{Introduction and problem statement} 
Given  $p$ symmetric  positive  definite matrices  $C_1, C_2,  \cdots,
C_p$,  each of  size $n\times  n$ the  standard Orthogonal 
Joint  Diagonalization (OJD) 
problem consists  of finding a  common $n\times n$ unitary  matrix $Q$
such that each $Q^T C_i Q$ is as close  as possible to a diagonal
$D_i$.
Note at the outset that in the `exact' case when there is a matrix 
$Q$ such that $C_i = Q D_i Q^T$  for $i=1,\cdots, p$, then the columns
of $Q$ are the eigenvectors of any one of the matrices, so in theory
the solution is completely determined by any one matrix in the set.
Issues of uniqueness have been studied in \cite{Asfari08}.
Here, uniqueness is meant up to permutations 
and sign scalings of the columns of $Q$.

An important source of 
joint  diagonalization problems  is in   Blind  Source
Separation  (BSS)  and independent  component  analysis  (ICA).
A typical  example can be described  as follows. Let
$\textbf{s}(t)$  be  an $M$-dimensional  vector  containing the  time
series of $M$ mutually independent source  signals that are mixed by a
mixing     matrix     $\textbf{A}^{n\times      M}$,    so    that
$\textbf{x}(t)=\textbf{A}\textbf{s}(t)$ is  an $n$-dimensional vector
containing the time series of $n$ mixture signals. In BSS one aims to
recover the source signals  $\textbf{s}(t)$ from the observed mixtures
$\textbf{x}(t)$ without  knowing $\textbf{A}$  or the  distribution of
the sources  but with the assumption
 that  the sources signals  are nonstationary
and mutually independent. The set  of matrices to joint diagonalize in
this  case  can be  obtained  from  covariance matrices  estimated  on
different time  windows or intercorrelation matrices  with time shifts
\cite{Belouchrani97}.  The  diagonal  matrices   in  this  case  correspond  to
covariance  matrices  or  intercorrelation   matrices  of  the  source
signals.  These matrices  are time  varying but  have always  diagonal
structures. Alternative approaches  to construct the set of
matrices  to  joint  diagonalize  can  be  obtained  by  using  linear
combinations of  intercorrelation matrices \cite{Vollgraf},  slices of
higher order cumulants tensors \cite{Cardoso}, time-frequency matrices
\cite{Fadaili}  or  second  derivatives   of  the  log  characteristic
function \cite{Yeredor}.

The various  algorithms that  have been
developed to solve this problem differ in the way in which the
actual optimization  problem is  solved and  what restrictions  on the
mixing matrix this  implies. 
They can be classified  into 
\emph{orthogonal}  joint diagonalization
  and  \emph{non-orthogonal}    joint diagonalization (NOJD) 
methods.  The most  popular  
algorithm for the OJD problem 
is  JADE \cite{Cardoso2},   a  Jacobi-like algorithm
based  on plane  rotations. 
In  BSS   applications,  the
assumption of  orthogonal mixing  and demixing  matrices can  often be
inappropriate.   In   such   situations,  methods   that   work   with
nonorthogonal        matrices       are        preferred
\cite{Yeredor2}\cite{Vollgraf2}\cite{Souloumiac}\cite{Aissa}. 

Most previous  work  in  this  area  focused  on  the  small  (dense)  case,
see, e.g., 
\cite{VdVeen01,CardosoSoul96,BGByersM93,FluryGautschi86}, among many
others.
In large applications where
 the signal vectors can consist of hundreds or more of
variables, existing joint  diagonalization algorithms
cannot be  applied due  to  their exorbitant  computational cost  or
simply because  they diverge for  high dimensional data sets. 
The case of large matrices has important applications in which
the $C_i$'s are often covariance matrices.
Thus, the article \cite{FTheisColor10} discusses `subspace coloring'
to tackle the autocorrelation structure
They propose a dimension reduction method that evaluates the
time structure of  multivariate observations, see also
\cite{Belouchrani97}.
  Their goal is to differentiate
the signal subspace from noise by extracting  a subspace of non
trivially autocorrelated data. In other words, their method is to find
a common invariant subspace for a set of correlation matrices $R_\tau$ where
the $\tau$'s correspond to a few time-lags.

The first question we address in this paper is to define a natural
extension of the OJD for large dimensional cases. 
For this it is important to think in terms
of subspaces. In this regard,
the  algorithms  we  propose are  neither  orthogonal   or
non-orthogonal but  approximate orthogonal.  
Indeed, in  practical high
dimensional  data application,  the investigator  is often  interested
only in  a subset of variables  that lie 
in a subspace of dimension  $k\ll n$. The underlying 
assumption is that the 
selected subspace concentrates most of the  information (variability)
of the data set.  The matrices are to be jointly diagonalized
only in this subspace.

\section{Approximate Joint Partial Diagonalization}
A  common  formulation of the standard OJD
 is to  seek  an $n \times n$ orthogonal matrix $Q$ that 
minimizes 
\eq{eq:obj0} 
f_0(Q) = \sum_{i=1}^p \| \Off(Q^T C_i Q) \|_F^2 ,
\en
where $\Off(X) $ denotes the 
matrix $X$ with its diagonal entries replaced by zeros.
In  situations where the matrices are large, methods based on 
minimizing the  above objective function over unitary matrices
$Q$,   become too expensive,  their costs being of the order
$O(p n^3)$ where $n$ is the size of the matrices. 
Our first goal is to 
define a generalization of the problem that can cope with large,
possibly sparse, matrices, while keeping the motivation and
rationale of the small dense case.

\subsection{Joint eigenspaces} 
The main change needed to  generalize
the joint diagonalization problem to large matrices
is to require that the orthogonal 
matrix $Q$  be of size $n \times k$ where $k$ is small.
In addition, \emph{we would also like $Q$ to represent the basis of  
  a common nearly invariant subspace for all $C_i$'s.}
With this, it is clear that a straightforward extension
of the formulation represented by \nref{eq:obj0} is not too meaningful.
Consider a naive generalization of the problem the problem that  would consist 
of  simply finding an orthogonal matrix 
(i.e., a matrix with orthonormal columns)
 $Q \in \RR^{n \times k}$ ($k$ columns)  so as to minimize
the same objective function \nref{eq:obj0}. This clearly does not work
in that it will not yield a basis of a subspace with the
desired property of near invariance 
by all matrices $C_i$. For example, consider the following 3 matrices
$C_i, i=1,2,3$ 
(only lower parts are shown due to symmetry) and assume we take $k=2$.
\[
\left[ 
\begin{array}{cc|ccc}
x   &      &    &   &    \\ 
    & x    &   &    &    \\ 
\hline 
x   &      & x &   &    \\ 
    & x    & x  & x  &   \\ 
x   &      &   & x & x 
\end{array}
\right]  \quad
\left[ 
\begin{array}{cc|ccc}
x   &      &    &   &    \\ 
    & x    &   &    &    \\ 
\hline 
x   &      & x &   &   \\ 
   x & x    &   & x  &   \\ 
    &  x    &   & x &  x 
\end{array}
\right]  
\quad
\left[ 
\begin{array}{cc|ccc}
x   &      &     &   &    \\ 
    & x    &     &   &    \\ 
\hline 
x   &       & x  &   &    \\ 
    &  x    &   & x  &    \\ 
x   &  x    & x  &  &  x 
\end{array}
\right]  
\]
In this example  all $C_i$'s happen to have
a $k \times k$ leading block $C_i(1:k,1:k)$ that is diagonal. Here, 
the optimal solution to the problem \nref{eq:obj0} is clearly the matrix 
$Q$ consisting of the first $k$ columns 
of the $n \times n$ identity matrix.
If entries in  the blocks (2,1) (and (1,2) by symmetry) were small in 
magnitude and if diagonal entries on the (1,1) block of each $C_i$
are near the largest eigenvalues of  $C_i$,
then  this $Q$ would be a good solution because  it would
represents the basis of a subspace that is near  the dominant invariant 
subspace for each of the 3 
matrices. Otherwise, i.e., in the case where entries in the $(1,2)$ blocks
are large, this optimal solution is not too meaningful.

Therefore, \emph{the primary problem we address in this paper is
to find  a matrix $Q$  that is of  size $n\times k$,  with orthonormal
columns such that the span of $Q$ is nearly invariant by each $C_i$, i.e.,
such that $C_i Q \approx Q D_i $ for some $k\times k$ matrix $D_i$.} 
In addition, we would like this 
common invariant subspace to be `dominant' for each $C_i$, i.e.,
associated with the largest eigenvalues for each $C_i$.
Note that the requirement that each $D_i$ be diagonal,
which was imposed in the standard case,  is being
omitted for now. Later, we will see that we can obtain the desired
near-diagonal structure in a post-processing
phase once the common invariant subspace has been identified.

One of several ways of formulating  the problem
rigorously, is to state  that we seek an orthogonal matrix 
$Q \in \RR^{n\times k}$ and $k \times k$ matrices $D_i$, for $i=1,\cdots,p$ 
that  minimize the objective function  
\eq{eq:obj} f(Q,D_1,...,D_p) = \sum_{i=1}^p \|
C_i  Q   -  Q  D_i  \|_F^2 .  \en
Consider each of the terms in the sum \nref{eq:obj}. A simple version
of Theorem 2.6 in \cite{Stewart-eig-book}  shows that once
$Q$ is known, each $D_i$ is uniquely determined. Indeed,
let $C \, \in \, \RR^{n\times n}$, \ $D\, \in \, \RR^{p\times p}$ 
and $R = C Q - Q D$ where  
$Q \in \, \RR^{n\times p}$  is  a (fixed)  orthogonal matrix. 
 Then the theorem states that 
in any unitarily invariant norm  $\| R \| $ is minimized
when $D = Q^T C Q$.
The  matrix $CQ - Q (Q^T C Q)$ is treated as a residual matrix 
when $Q$ is considered as an approximate eigenspace. This theorem
will play a key role in the analysis and so we show a slightly more elaborate
version in Section~\ref{sec:genres}.

Since each of the norms $\| C_i Q - Q D_i \|_F$ represents a sort of 
measure of the invariance of the span of $Q$ by $C_i$, minimizing
\nref{eq:obj} will produce a  common near eigenspace  for  the  $C_i$'s.  
In realistic applications, such  a subspace does not exist,  but there is
a  common subspace  of  small  dimension  $k$ that  is  `almost
invariant' by  all $C_i$'s may arise  naturally in some applications.
 This  translates into the existence  of a
matrix $Q$  (orthogonal basis  of this  subspace) and  matrices $D_i$,
$i=1,\cdots,p$ such  that each for  the norms $\|C_i  Q - Q  D_i\|$ is
small.

\subsection{Joint partial diagonalization} 
The requirement that the $D_i$'s  be diagonal has been  relaxed so far
but we can make the $D_i$'s
diagonal in a \emph{post-processing phase}
by a `rotation' of the bases. 
In a trivial  approach to  the problem, we could  just diagonalize
each $D_i$ separately as $D_i = W_i \Lambda_i W_i^T$,
where $W_i \in \RR^{k \times k}$, $W_i^T W_i = I$, 
and $\Lambda_i$ is diagonal. Then
\eq{eq:invar1}
\| C_i Q - Q D_i \|_F = \|  C_i (Q W_i) -(Q W_i)  \Lambda_i \|_F .
\en
Thus, each $C_i$ has been partially diagonalized with 
orthogonal matrices of the form $Q W_i$, 
where each $W_i$ is $k\times k$ and unitary.
These orthogonal bases $Q W_i$ span the same space which is a
nearly invariant subspace for each $C_i$ as desired but each
orthogonal basis $Q_i\equiv Q W_i$ 
is different for each  $C_i$'. This
 solution answers  a slightly different question
from the standard one. 

It is possible to make all matrices $Q ^T C_i Q$ nearly diagonal with the same
matrix $Q$ as is done in the classical case, to obtain in this way an
approximate common partial diagonalization.  A first approach to the proble is
to rely on orthogonal transformations.  Let $D_i = Q ^T C_i Q$ for
$i=1,\cdots,p$ and assume that \emph{the $D_i$'s can be approximately jointly
  diagonalized}. Then there exists a matrix $W\in \RR^{k \times k}$, with
$W^T W = I$, such that $D_i W \approx W \Lambda_i $ for $i=1, 2, \cdots, p$,
where $\Lambda_i$ is diagonal. If, in addition, $C_i Q \approx Q D_i $, then
\[ C_i Q W \approx (Q W) W^T D_i W \approx (Q W) \Lambda_i . \] This results in
an approximate joint diagonalization of the $C_i$'s by the orthogonal matrix
$Q W$. We will later analyze the error made in this 2-step process and show that
the square of the Frobenius norm of the error is just the sum of the square of
the Frobenius norms of the errors made in the two steps.

A second approach, which can be termed non-orthogonal, is seek a  partial
diagonalization of
the $C_i$'s by congruence transformations. From the first step, we
obtain the approximation:
\[ Q^T C_i Q \approx D_i . \]
 Then in a second step, a common 
$k\times k$ congruence transformation is applied to the $D_i$'s so
that $ L^T D_i L \approx \Lambda_i$ and therefore:
\[
L^T Q^T C_i Q L \approx \Lambda_i . \]
Now we have a basis, namely $QL$,
 of a nearly invariant subspace in which each $C_i$ has 
a nearly diagonal representation.
We will not consider non-orthogonal transformations in this paper
but wish to emphasize here that it is also possible to 
partly diagonalize the $C_i$'s with non-orthogonal transformations.

\subsection{General approach} 
The main method introduced in this paper solves
the decoupled formulation of the problem as described above.
In a first stage, the objective function \nref{eq:obj} is minimized 
with respect to $Q$ and the $D_i$'s only, ignoring the diagonal structure 
requirement for the $D_i$'s. This yields an orthogonal matrix $Q$
and matrices $D_i$, 
so that 
$C_i Q \approx Q D_i$, for $i=1,\cdots, p$, i.e., $Q$ represents
 the basis of a nearly invariant subspace common to all
the $C_i$'s. 

If desired, the resulting approximate factorizations
$C_i Q \approx Q D_i$, can be post-processed in a second stage
 to get a joint partial diagonalization of the $C_i$'s.
This is done via the solution of a standard,
orthogonal or non-orthogonal,  joint diagonalization of
the $D_i$'s as described in the previous subsection. 
This  standard joint diagonalization problem  involves
the matrices $D_i$, which are 
of size $k \times k$ where $k$ is typically small.

\section{General results}\label{sec:genres} 
The objective function \nref{eq:obj}
is not convex with respect to $Q, \{D_i\}$ together.
It is however quadratic with respect to $Q$ alone, and each $D_i$ alone.
As it turns out we can optimize the objective by focussing on $Q$ alone.
The $D_i$'s will be obtained immediately once the optimal $Q$ is
known. This section discusses this and related issues.

Recall that when we consider matrices as vectors in $\RR^{n^2}$, 
the usual Euclidean inner product translates into the Frobenius
inner product of matrices
\eq{eq:FroDP}
\left\langle A, B \right\rangle = \trace \ [B^T A].  
\en 
In particular $\| X \|_F^2 = \trace [X^T X] $. 
In addition, we will say that $A \perp_F B $ if 
$\left\langle A, B \right\rangle = 0$.
If $Q$ is an $n\times k$ orthogonal matrix then the orthogonal projector
onto the span of $Q$, which we will call $P$,
 is represented by the matrix $P=Q Q^T$. 
Hidden in this notation is the fact that the same $P$ 
is represented in this form for different orthogonal bases of the same subspace.

\subsection{Problem decoupling} 
The first question we address is: 
Assuming that $Q$ is fixed, what is the best set of 
$D_i$'s that minimize \nref{eq:obj}? 
Note that the resulting optimization problem decouples into 
$p$ subproblems. Indeed,   in order to minimize 
\[  \sum_{i=1}^p \| C_i Q - Q D_i \|_F^2 
\]
with respect to the $D_i$'s, 
we can minimize each of the terms separately. The
following result is a reformulation of 
Theorem~2.6 in \cite{Stewart-eig-book}.

\begin{theorem} 
For a fixed orthogonal matrix $Q$, define $D_{Q,i} = Q^T C_i Q$.
Then the objective function \nref{eq:obj}
is minimized when $D_i = D_{Q,i}$. In addition, if 
we define the residual matrix as 
$R_i = C_i Q - Q D_{Q,i} $ then 
$R_i = (I-Q Q^T) C_i Q$ and in particular
$R_i \perp_F Q$ and 
$\|R_i \|_F = \| (I- Q Q^T) C_i Q \|_F$.
\end{theorem}
\begin{proof}
We   start with  the following relation:
\begin{eqnarray*}
 C_i Q - Q D_i  &=&  Q Q^T C_i Q  + (I-QQ^T) C_i Q - Q D_i \\
  &=&    Q (Q^T C_i Q - D_i)        +    (I-QQ^T) C_i Q .
\end{eqnarray*}
Exploiting the orthogonality of matrices of the form $Q X $ and $(I-Q Q^T)Y$
 we obtain:
\[ 
\| C_i Q - Q D_i \|_F^2  =
    \| Q (Q^T C_i Q - D_i) \|_F^2       +    \|(I-QQ^T) C_i Q \|_F^2
\]
which  shows  that the minimum is indeed achieved for $D_i = Q^T C_i Q$.
In this case $R_i =(I-QQ^T) C_i Q$. As a result it is clear that
$R_i \perp_F Q$ 
and that the error term relative to $C_i$ is 
$ 
\| R_i \|_F = 
\| (I-QQ^T) C_i Q\|_F. 
$ 
This completes the proof.
\end{proof} \\
We emphasize an important result of the theorem which 
is that the residual 
$R_i$ is orthogonal
to $Q$ in the Frobenius inner product sense, i.e.,  $\trace (Q^T R_i) = 0$.  
However, it is clear that the stronger result 
$Q^T R_i = 0$ is also true.

As mentioned in the introduction, we can first 
optimize \nref{eq:obj} with respect to $Q$ only and then find an optimal joint
diagonalization for the $D_i$. Both steps incur an error and it
turns out that these two errors decouple. 

\begin{proposition}[Error decoupling]
Assume that we find a joint orthogonal matrix $Q$ such that 
\eq{eq:erS}
C_i Q - Q D_i = E_i \quad i=1,\cdots, p 
\en
where $D_i = Q_i^T C_i Q_i$, 
and  that an approximate joint diagonalization 
of the matrices $D_i$'s is available with:
\eq{eq:erK}
D_i W - W \Lambda_i = F_i , \ \quad i=1,\cdots, p .
\en
Define $\tilde Q = Q W$.
Then, 
\eq{eq:erT}
\sum_{k=1}^p \| C_i \tilde Q - \tilde Q \Lambda_i \|_F^2 
= 
\sum_{k=1}^p\left[ \| E_i \|_F^2 + \| F_i \|_F^2 \right] .
\en
\end{proposition}
\begin{proof}
For each $i$:
$ C_i \tilde Q - \tilde Q \Lambda_i = C_i Q W -  Q W \Lambda_i 
= C_i Q W -  Q (D_i W -F_i)$. 
Hence 
\begin{eqnarray*}
\sum_{i=1}^p \| C_i \tilde Q - \tilde Q \Lambda_i \|_F^2 
&=& 
\sum_{i=1}^p \| (C_i Q W -  Q D_i W ) + Q F_i \|_F^2 
\\
&=& 
\sum_{i=1}^p \|[(C_i Q - Q D_i)  + Q F_i W^T ] W \|_F^2
\\
&=& 
\sum_{i=1}^p \| E_i + Q F_i W^T \|_F^2 \ . 
\end{eqnarray*}
Since each
$E_i = C_i Q - Q D_i $ is orthogonal to $Q $ in the sense of
the Frobenius inner product we have: 
\begin{eqnarray*}
\sum_{i=1}^p \| E_i + Q F_i W^T \|_F^2
&=& \sum_{i=1}^p \|E_i \|_F^2   + \| Q F_i W^T  \|_F^2 \\
&=&   \sum_{i=1}^p \| E_i \|_F^2   + \|  F_i \|_F^2 .
\end{eqnarray*}
\end{proof}

If the errors decouple, can we say that the two optimization problems
 decouple as well? As it turns out the matrices $\Lambda_i$ for the second
optimization problem (JOD of the $D_i$'s) 
are also determined from the knowledge of
the optimal $Q$. 
\begin{theorem}\label{th:equiv0} 
Let $\ubarQ$ an $n\times k$ orthogonal matrix and 
$\ubarLambda_i$,  $i=1,\cdots,p$ diagonal matrices such that
$f(\ubarQ,\ubarLambda_1, \cdots, \ubarLambda_p)$ is minimum. 
Then $\ubarLambda_i = \diag(\ubarQ^T C_i \ubarQ) $ for $i=1,\cdots, p$.
In addition, for this $\ubarQ$ the objective function \nref{eq:obj} is equal to:
\eq{eq:obj2}
f(\ubarQ,\ubarLambda_1, \cdots, \ubarLambda_p) =
g(\ubarQ) \equiv   \sum_{i=1}^p 
\left[  \| C_i \ubarQ - \ubarQ D_{\ubarQ,i}\|_F^2 +  \|\Off( D_{\ubarQ,i}) \|_F^2 \right] , 
\en
where $D_{\ubarQ,i} = \ubarQ^T C_i \ubarQ$.
\end{theorem} 
\begin{proof} 
Consider one term in the sum \nref{eq:obj} where $\ubarQ$ is optimal.
Dropping the index $i$ we denote by $D_{\ubarQ}$ the matrix $\ubarQ^TC\ubarQ$ and recall
that $\ubarQ^T(C\ubarQ-\ubarQ D_\ubarQ) = 0$. Then,
\begin{eqnarray*}
\| C  \ubarQ -  \ubarQ \ubarLambda \|_F^2 
&=& 
\| C \ubarQ -  \ubarQ D_{\ubarQ} + \ubarQ ( D_{\ubarQ} -  \ubarLambda) \|_F^2\\
&=& 
\| C \ubarQ -  \ubarQ D_{\ubarQ} \|_F^2 + 
\| \ubarQ (D_{\ubarQ} - \ubarLambda) \|_F^2 \\
&=& 
\| C \ubarQ -  \ubarQ D_{\ubarQ} \|_F^2 + 
\| D_{\ubarQ} - \ubarLambda \|_F^2 .
\end{eqnarray*}
The term $\|D_{\ubarQ} - \ubarLambda \|_F^2$ is smallest when $\ubarLambda = \diag (D_{\ubarQ})$.
This shows the first part. The second part follows by observing that
 for the optimal $\ubarLambda$ we have 
 $\|D_{\ubarQ} - \ubarLambda \|_F^2 = \| \Off (D_{\ubarQ}) \|_F^2$.
\end{proof} 

The theorem shows that the problem of  minimizing 
$f(Q,\Lambda_1, \cdots, \Lambda_p)$ where the $\Lambda_i$'s are diagonal,
can be restated entirely in terms of the unknown matrix $Q$:
it suffices to 
\emph{find an orthogonal $Q$ so that \nref{eq:obj2}
  is minimum. } Each $\Lambda_i$ is then equal to $\Diag(D_{Q,i})$.

The objective function $g$ defined in \nref{eq:obj2} 
  can be split in two parts:
\eq{eq:obj3}
\gamma(Q) =  \sum_{i=1}^p \|C_i Q - Q D_{Q,i}\|_F^2; \qquad 
\omega(Q) =  \sum_{i=1}^p  \|\Off( D_{Q,i}) \|_F^2 .
\en
Now we can think of a 2-stage, or decoupled,  procedure.
The first stage of this procedure finds 
an optimal $ Q$ that  minimizes $\gamma(Q)$. The second stage
seeks a basis change $Q_{new} = Q W$, where the desired
$W \in \RR^{p\times p} $ is unitary, and is selected so as  to minimize
$\omega(Q_{new})$. 
As the following Lemma shows, this second stage does not change
the minimal value of $\gamma$ obtained in the
first stage.
\begin{lemma}\label{lem:qw}
Let $W$  be any  $k\times k$ unitary matrix and let
$\gamma(Q) $ defined in \nref{eq:obj3}. Then
\eq{eq:invar2} 
\gamma(Q) = \gamma(QW) . 
\en
\end{lemma} 
\begin{proof} 
We have
\begin{eqnarray*}
\|C_i Q - Q D_{Q,i}\|_F^2 &=& \| C_i Q W - Q D_{Q,i}W \|_F^2 \\
&=& \|C_i (QW) - (QW) W^T D_{Q,i}W \|_F^2\\
&=& \|C_i (QW) - (QW) [(QW)^T C_i (QW)]  \|_F^2 .
\end{eqnarray*}
This completes the proof by noting that 
$(QW)^T C_i (QW) \equiv D_{QW,i} $.
\end{proof}\\
The lemma  establishes that $\gamma(Q)$ 
is a function of the subspace only, not its basis $Q$. 

Consider now a certain $\tilde Q$ found in a first stage to minimize
$\gamma(Q)$. In the second stage  we wish to minimize 
$\omega(Q) $ over different orthonormal bases  $Q= \tilde Q W$ of the 
same subspace.  Note that for any unitary $W$, we have
\[
D_{QW,i} = (QW)^T C_i (QW)  = W^T D_{Q,i} W . \]
As a result, the best  $W$ can be found by  jointly minimizing
 the second term of the objective function, namely:
\eq{eq:obj4}
\sum_{i=1}^p \|\Off [ W^T D_{\tilde Q,i} W ] \|_F^2,  
\en 
over $W$. This is the classical joint diagonalization problem 
but it now involves small matrices  of size $k\times k$.

\begin{proposition}
A minimizer $Q_*$ of the objective function \nref{eq:obj2}
is of the form $Q_* = \tilde Q W$ where $\tilde Q$ minimizes the function
$\gamma(Q)$ over all orthogonal $m\times p$ matrices, and 
$W$, minimizes  \nref{eq:obj4}, i.e., it solves a standard 
orthogonal joint diagonalization problem for the matrices
$D_{\tilde Q,i} = \tilde Q^T C_i \tilde Q_i$, $ i=1,\cdots,p$.
\end{proposition} 

\begin{proof} 
The proof follows from the lemma and the discussion above.
\end{proof} \\
The 2-stage procedure just described is sketched below as
Algorithm~\ref{alg:decoupl}.
 
\begin{algorithm}
  \caption{Decoupling procedure}\label{alg:decoupl} 
  \begin{algorithmic}[1]
\State   Find an orthogonal $Q$ that minimizes 
$\gamma(Q) = \sum_{i=1}^p \|C_i Q - Q D_{Q,i}\|_F^2 $  

\State   Find a unitary matrix $W \ \in \RR^{k\times k}$ 
that minimizes  \nref{eq:obj4}  

\State   Output $\tilde Q \equiv QW $ and $\tilde \Lambda_i = \diag [
\tilde Q^T C_i \tilde Q ] 
= \diag [ W^T D_{Q,i} W ]$ 
\end{algorithmic}
\end{algorithm}

Observe that the algorithm will yield an orthogonal matrix  $\tilde Q$ and
diagonal matrices $\tilde \Lambda_1,\cdots, \tilde \Lambda_p$ such that
$\gamma (\tilde Q) = \gamma (QW)$ is minimum over orthogonal matrices $Q$
and, exploiting orthogonality as in  the proof of Theorem~\ref{th:equiv0} we get
\begin{eqnarray*}
\sum_{i=1}^p \|C_i \tilde Q - \tilde Q \tilde \Lambda_i \|_F^2
  &=&\sum_{i=1}^p \| [C_i \tilde Q - \tilde Q D_{\tilde Q_i} ] + 
            \tilde Q [D_{\tilde Q_i} -\tilde \Lambda_i]\|_F^2 \\
  &=&\sum_{i=1}^p \|C_i \tilde Q - \tilde Q D_{\tilde Q_i}\|_F^2
      + \sum_{i=1}^p \| \tilde Q [D_{\tilde Q_i} -\tilde \Lambda_i]\|_F^2\\
  &=& \gamma(\tilde Q) + \sum_{i=1}^p \| D_{\tilde Q_i} -\tilde \Lambda_i\|_F^2
  =  \gamma(\tilde Q) + \sum_{i=1}^p \| \Off [D_{\tilde Q_i} ]\|_F^2
  \end{eqnarray*}
The last equality comes from the definition of $\tilde \Lambda_i$.
  This establishes the relation
  \eq{eq:g1tq}
f(\tilde Q, \{ \tilde \Lambda_i\}_{i=1:p} ) = 
  \sum_{i=1}^p \|C_i \tilde Q - \tilde Q \tilde \Lambda_i \|_F^2
  = \gamma(\tilde Q) +  \omega(\tilde Q)
  \en
which shows a similar result to that of Theorem~\ref{th:equiv0} for $\tilde Q$.

We may now ask whether or not  the above procedure can
recover a solution that is the same or
close to the optimal one for \nref{eq:obj} in which 
the additional constraint that each $D_i$ be diagonal is enforced. 
Let $\tilde Q$ be the matrix   that results from
the decoupling procedure of Algorithm~\ref{alg:decoupl}.
At the same time let $\ubarQ, \ubarLambda_{1}, \cdots, \ubarLambda_{p}$
be a minimizing set of matrices for  \nref{eq:obj} as defined in
Theorem~\ref{th:equiv0}. Then
because $\tilde Q $ minimizes $\gamma$  we  clearly have
\eq{eq:ineq0}
\gamma(\tilde Q) \le \gamma (\ubarQ) \ .
\en

Now, from Theorem~\ref{th:equiv0}, the  optimal value of the 
 objective function is equal to:
 \eq{eq:ineq1}
 f(\ubarQ, \ubarLambda_{1}, \cdots, \ubarLambda_{p}) = \gamma(\ubarQ) + \omega(\ubarQ)
 \en
 On the other hand, the set $\tilde Q, \{ \tilde \Lambda_i \}_{i=1:p} $
 realizes a  partial orthogonal diagonalization and so:
 \eq{eq:ineq2}
 f(\ubarQ, \ubarLambda_{1}, \cdots, \ubarLambda_{p}) \le
 f(\tilde Q, \tilde \Lambda_{1}, \cdots, \tilde \Lambda_{p})
 \en
 Putting the above relations together
 will give the following corollary.

 \begin{corollary}
   With   the  notation of  Theorem~\ref{th:equiv0} and
   Algorithm~\ref{alg:decoupl}, the following double inequality is satisfied:
   \eq{eq:ineq3}
   f(\tilde Q, \{\tilde \Lambda_i\}_{i=1:p} )
   - [\omega(\tilde Q) -\omega(\ubarQ)]
 \le f(\ubarQ, \{\ubarLambda_i\}_{i=1:p} )
 \le f(\tilde Q, \{\tilde \Lambda_i\}_{i=1:p} )  \ . 
\en
\end{corollary}
\begin{proof}
  The right part of \nref{eq:ineq3} is just \nref{eq:ineq2}.
  From \nref{eq:ineq1} and \nref{eq:ineq0} we get
 \[
 f(\ubarQ, \ubarLambda_{1}, \cdots, \ubarLambda_{p}) = \gamma(\ubarQ) + \omega(\ubarQ)
 \ge \gamma(\tilde Q) + \omega(\ubarQ) .
\]
This yields the first inequality after substituting  the value of
$\gamma(\tilde Q) $  extracted from equality \nref{eq:g1tq}.
  \end{proof}


  A by-product of the result is that the difference $\omega(\tilde Q)-\omega(\ubarQ)$
  is nonnegative.  It is zero when the optimal solution is also the optimal
  solution obtained by the decoupling algorithm.  This happens when there exists
  an \emph{exact} joint diagonalization of the projected set
  $\{ D_{\tilde Q, i} \}_{i=1,\cdots,p}$. In this situation,
  $\omega(\tilde Q) = 0$ which forces $\omega(\ubarQ)$ to also equal zero since
  $\omega(\tilde Q) \ge \omega(\ubarQ)$ and then \nref{eq:ineq3} implies that
  $ f(\ubarQ, \{\ubarLambda_i\}_{i=1:p} ) = f(\tilde Q, \{\tilde
  \Lambda_i\}_{i=1:p} )$.  In case the matrices $D_{\tilde Q,i}$ can be jointly
  diagonalized only inexactly, then the optimal objective functions obtained by
  the exact optimum and the decoupled one differ by $\omega(\tilde
  Q)-\omega(\ubarQ)$. This difference will be small if $\omega(\tilde Q)$ is small,
  i.e., if the set $ \{ D_{\tilde Q, i} \}$ is nearly jointly diagonalizable.
 
  As a final note, we should mention that no statement was made on the closeness
  of the subspace obtained by the decoupling procedure to the optimal one.  The
  lack of uniqueness implies that even if the related objective functions are
  equal, we cannot say that the subspaces are the same or close.

\subsection{The Grassmannian perspective}\label{sec:grass} 
As  already stated, Algorithm~\ref{alg:decoupl} in the
simple form given above does not specify which of the many possible
nearly invariant subspaces is selected when minimizing $\gamma(Q)$.
In typical applications, it is the  subspace 
associated with the largest $p$ eigenvalues for each $C_i$ that is desired.
We could consider instead the objective function:
\eq{eq:Phi} \phi(Y) = \half \sum_{i=1}^p \trace [ Y^T C_i Y ] .
\en
We changed notation from $Q$ to $Y$ to reflect a  more generally adopted notation in this context. 
In the case when $p=1$, it is well-known \cite[10.6.5]{GVL-book}
that maximizing the above trace over orthogonal matrices $Y$ will yield an
orthogonal basis for the invariant subspace associated with the largest $p$
eigenvalues. This is exploited in the TraceMin
algorithm~\cite{SamehTong00,Sameh-trace-min} a  method designed 
 for computing an invariant subspace associated with smallest 
eigenvalues for standard and generalized eigenvalue problems. 

Unfortunately when $p>1$ the optimum is  just the dominant subspace
of the matrix $\sum_{i=1}^p C_i$. Indeed, maximizing \nref{eq:Phi}
is oblivious to the 
invariance of $[Y]$ with respect to each matrix $C_i$. This is just the
opposite of what is obtained when minimiming $\gamma(Y)$:
optimize invariance regardless of which subspace is obtained.
What is needed is
a criterion that mixes the two objective functions.
An objective function  similar to \nref{eq:Phi} was also proposed  in
\cite{FTheisColor10} where the trace in \nref{eq:Phi} is replaced by
the Frobenius norm squared:
\[
  \psi(Y) =  \half \ \sum_{i=1}^p \| Y^T C_i Y \|_F^2 .
\]
Let us  examine this alternative. First note that due to symmetry,
$ \| Y^T C_i Y \|_F^2  = \half \trace [ Y^T C_i Y  Y^T C_i Y ] $. Then letting $\Pi = Y Y^T$,
we have:
\begin{align*}
  \trace (  Y^T C_i \Pi  C_i Y  ) 
  &=   \trace (  Y^T C_i (I- (I-\Pi))  C_i Y  ) \\
  &=   \trace (  Y^T C_i^2 Y  ) - \trace (  Y^T C_i (I-\Pi)  C_i Y  ) \\
  &=   \trace (  Y^T C_i^2 Y  ) - \trace (  Y^T C_i (I-\Pi) (I-\Pi)  C_i Y ) \\
  &=   \trace (  Y^T C_i^2 Y  ) - \| (I-\Pi)  C_i Y  \|_F^2 
\end{align*}
and upon summation we therefore get: 
\eq{eq:Psi0}
\psi(Y) =  \half \sum_{i=1}^p \trace (  Y^T C_i^2 Y  )  - \half \gamma(Y) . 
\en
Thus, this alternative provides the desired mixing of the objective functions: by
maximizing $\psi(Y)$ we maximize the sum of the traces of $Y^T C_i^2 Y$ while at the same
time making the invariance mesure $\gamma(Y)$ small as desired.
Note that we could also replace the matrices $C_i^2$ in the first term by $C_i$ as an
alternative and this would change the first term to $\phi(Y)$ defined in \nref{eq:Phi}, so
$\psi$ could be defined as $\psi(Y) = \phi(Y) - \gamma(Y)$.

To  better balance the criteria of invariance (small $\gamma$) and subspace targetting
(large $\phi$),  we   introduce a parameter $\eta $ and redefine $\psi$ as:
\eq{eq:Psi}
\psi(Y) =   \phi (Y)  - \eta \  \gamma(Y) .
\en


A key observation in the definition~\nref{eq:Psi} is that $\psi(Y) $ is
invariant upon unitary transformations. In other words if $W$ is a $p \times p$
unitary matrix then, $\psi( Y W) = \psi(Y)$. This suggests that it is possible,
and possibly advantageous, to seek the optimum solution in the 
Grassmann manifold~\cite{eas:99}.
Recall, from e.g., \cite{eas:99}, that the Stiefel manifold is
the set of orthogonal transformations\footnote{This is often termed the \emph{compact}
  Stiefel manifold. The standard   Stiefel manifold has no orthogonality requirement for the
  columns of $Y$ but $Y$ must be of full rank.}:
\eq{eq:Stnp}
St(n,p) = \{ Y \ \in \ \RR^{n \times p} \ : \  Y^T Y = I \} . 
\en
while the Grassmann manifold is the quotient manifold
\eq{eq:Grnp}
G(n,p) = S(n,p) / O(p) 
\en
where $O(p)$ is the orthogonal group of unitary $p \times p$ matrices.
Each point on the manifold, one of the equivalence classes in the above definition,
can be viewed as a subspace of dimension $p$ of $\RR^n$. 
It can be indirectly represented by a basis $V  \in St(n,p)$
modulo a unitary transformation and so we will denote it by $[V]$, keeping in
mind that it does not matter which member $V$ of the equivalence class is selected
for this representation.

For a given $Y \in \ St(n,p)$ the tangent space of the Grassmann manifold at
$[Y]$ is the set of matrices $\Delta \in \RR^{n\times p}$ satisfying the
relation
\eq{eq:Tang}
Y^T \Delta = 0,
\en
see, \cite{eas:99}.


It is possible to  adapt the treatment in \cite{eas:99}  to our slightly different context
and obtain a Newton-type procedure on the Grassmann manifold.
However, a  drawback of this approach is that  it only seeks a
stationary  point on the manifold, one for which the  gradient vanishes.
The limit can be  any invariant subspace not the desired one and so this method will again
miss the objective of targetting a specific  subspace.

We will instead explore a Gradient based procedure.
Exploiting results in \cite{eas:99}, it is easy to see that
when expressed on the Grassmann manifold, the gradient of \nref{eq:Phi} at $[Y]$ is 
\eq{eq:gradPhi}
\nabla \phi_Y  =   \sum_{i=1}^p  (I-Y Y^T) C_i Y \equiv \sum_{i=1}^p  (C_i Y -Y D_{Y,i}) . 
\en
Similarly we now need the gradient of $\gamma(Y)$.
We write
\eq{eq:gamG}
\gamma(Y) =  \sum_{i=1}^p \gamma_i(Y), \quad \mbox{with} \quad
\gamma_i(Y) = 
\half \| R_i(Y) \|_F^2 ,  \quad R_i(Y) \equiv C_i Y - Y D_{i,Y} .
\en
We will look for the gradient of each $\gamma_i$.
Note in passing  that 
\eq{eq:Fi}
\nabla \phi_Y  =   \sum_{i=1}^p R_i(Y); \quad \mbox{and} \quad 
R_i(Y) = (I - Y Y^T) C_i Y   = (I-\Pi)C_i Y .
\en

\begin{proposition}
  The gradient of the objective function \nref{eq:Psi} at point $[Y]$ of
  the Grassmann manifold
  is given by~:
  \eq{eq:gradPsi}
  \nabla \psi_Y 
          \ = \  \sum_{i=1}^p [(1-\eta ) I + \eta \Pi  C_i] R_i(Y) +
          \eta \sum_{i=1}^p R_i(Y) D_{i,Y} . 
          \en
  \end{proposition} 

\begin{proof}
  For any $\Delta \ \in \RR^{n\times p}$ we consider
\begin{eqnarray*} 
  R_i(Y+\Delta)
  &=& (I - (Y+\Delta)  (Y+\Delta)^T ) C_i \ (Y+\Delta) \\
  &=& [(I - Y Y^T) - Y \Delta^T   - \Delta Y ^T  - \Delta\Delta^T Y ]
      C_i \ (Y+\Delta) \\
  &=& R_i (Y)  + (I - Y Y^T) C_i \Delta - (Y \Delta^T + \Delta Y^T ) C_i Y + O(\Delta^2) .
\end{eqnarray*}
Therefore,
\begin{multline}
  \half \|     R_i(Y+\Delta) \|_F^2 =
  \half \|     R_i(Y) \|_F^2 + \\
  \left\langle R_i(Y),  (I - Y Y^T) C_i \Delta - (Y \Delta^T + \Delta Y^T ) C_i Y \right\rangle + O(\Delta^2)
 \end{multline}
 Noting that $R_i(Y) $ is orthogonal to $Y$, this becomes:
 \[
   \gamma_i (Y+\Delta) = \gamma_i (Y) +
     \left\langle R_i(Y),  (I - Y Y^T) C_i \Delta - \Delta Y^T C_i Y \right\rangle + O(\Delta^2)
   \]
   Thus, noting that $<A, B\Delta> = <B^T A, \Delta>$ and
   $<A, \Delta B> = < A B^T, \Delta > $ : 
   
   \begin{eqnarray*}
     <\nabla \gamma_i, \Delta>
     &=&  \left\langle R_i(Y),  (I - Y Y^T) C_i \Delta - \Delta Y^T C_i Y \right\rangle \\
     &=&  \left\langle C_i (I-\Pi) R_i(Y),  \Delta \right\rangle
     - \left\langle R_i(Y),  \Delta D_{i,Y} \right\rangle \\
     &=&  \left\langle C_i R_i(Y),  \Delta \right\rangle
     - \left\langle R_i(Y) D_{i,Y},  \Delta \right\rangle \\
     \end{eqnarray*} 
     and therefore since this equality must be true from any $\Delta $ on
     tangent space $\Delta = (I-\Pi)\Delta$ the gradient of $\Gamma_i$ becomes
     $
     \nabla \gamma_i = (I-\Pi) C_i R_i(Y) - R_i(Y) D_{i,Y} $ from which we get
     \eq{eq:DelGam}
     \nabla \gamma_Y =
     \sum_{i=1}^p \left[ (I-\Pi) C_i R_i(Y) - R_i(Y) D_{i,Y} \right]
     \en
     and, using \nref{eq:Fi} the expression of the gradient of $\psi$ is:
      \begin{eqnarray*}
        \nabla \psi_Y 
        &=&   \sum_{i=1}^p  R_i(Y) - \eta  \sum_{i=1}^p
            [(I-\Pi) C_i R_i(Y) - R_i(Y) D_{i,Y}] \nonumber \\
        &=&   \sum_{i=1}^p  \left[ R_i(Y)  - \eta (I-\Pi) C_i R_i(Y) +
            \eta  R_i(Y) D_{i,Y} \right] \label{eq:DelPsi} \\
        &=&   \sum_{i=1}^p [(1-\eta ) I + \eta \Pi  C_i] R_i(Y) +
     \eta \sum_{i=1}^p R_i(Y) D_{i,Y}  . 
      \end{eqnarray*}
      This is the desired expression.
    \end{proof}

      This expression will be exploited in Section \ref{sec:grad} to develop a gradient
      ascent type algorithm. 

  
\subsection{Comparison with the global optimum}
We now return to the `global' (full)
Joint Orthogonal Diagonalization problem with a goal 
of unraveling relationships between different approaches.
Here, we seek a unitary matrix $Q_g$ (now a member of $\RR^{n \times n}$) such
that $Q_g^T C_i Q_g $ is as close as possible to a diagonal matrix for all
$i$'s. The subscript $g$ helps to distinguish between orthogonal matrices that are
in $\RR^{n \times n}$ (denoted by $Q_g$) and those 
in $\RR^{n \times k}$ (denoted by $Q$).
There are two possible ways of formulating the global optimization 
problem. The first is  represented by the objective function 
\nref{eq:obj0} which we write for convenience as
\eq{eq:obj0D}
f_0(Q_g) = \sum_{i=1}^p \| Q_g^T C_i Q_g  - \Diag [ Q_g^T C_i Q_g  ] \  \|_F^2 . 
\en
Since $Q_g$ is unitary, we can rewrite $f_0$ as:
\eq{eq:obj03}
f_0(Q_g) = \sum_{i=1}^p \| C_i Q_g  -  Q_g \Diag [ Q_g^T C_i Q_g  ] \  \|_F^2 . 
\en
The second formulation is simply \nref{eq:obj} in which $Q$ is
now a $n\times n$ matrix $Q_g$ and we add 
the  constraint that the $D_i$'s be diagonal for $i=1,\cdots,p$.
A corollary of Theorem\ref{th:equiv0} is that 
 the two formulations yield identical solutions.
\begin{corollary} 
Consider the problem of minimizing \nref{eq:obj} 
rewritten as
\eq{eq:objG} f(Q_g,\Lambda_1,...,\Lambda_p) = \sum_{i=1}^p \|
C_i  Q_g   -  Q_g  \Lambda_i  \|_F^2 .  \en
under the additional
constraint that the $\Lambda_i$'s are all diagonal matrices.
Then, when the optimum is reached, the $\Lambda_i$'s must satisfy: 
\[ 
\Lambda_i = \Diag [ Q_g^T C_i Q_g ] \quad \mbox{for} \quad i=1,\cdots, p.
\]
Therefore, the three objective functions 
 \nref{eq:obj0D},  \nref{eq:obj03}, and \nref{eq:objG} have the same
minimum value and produce essentially the same solution. 
\end{corollary}
\begin{proof} 
  The first part is a simple consequence of Theorem~\ref{th:equiv0} which 
is clearly valid when $Q_g \in \RR^{n\times n}$.
It was shown above that \nref{eq:obj0D} and \nref{eq:obj03} are equivalent.
In addition, the first part of this corollary shows  that if the optimum of
\nref{eq:objG} is reached 
for  $Q_g, \{ \Lambda_i\}$, we must have $\Lambda_i = \Diag [Q_g^T C_i Q_g]$.
Therefore the minimum value of this objective function is the same as that 
of \nref{eq:obj03}.
\end{proof} \\
The term `essentially' in the statement of the result
refers to the fact that the minimizer is unique only up to  signs
for the columns of $Q_g$.

Assume now that we solve the problem by using \nref{eq:obj03} - 
or, equivalently, \nref{eq:objG}.
A rank $k$ is selected and the optimal matrix $Q_g$ is split as 
 $Q_g = [Q_1, Q_2]$ where $Q_1$ has $k$ columns. Each diagonal matrix
$\Lambda_i$ is also split accordingly into a $k \times k$ matrix $\Lambda_i\up{1}$,
associated with the dominant eigenvalues,  and
an $(n-k)\times (n-k)$ matrix $\Lambda_i\up{2}$ associated with the
remaining eigenvalues. 

Consider again one term in the sum \nref{eq:objG} which we split 
appropriately as follows:
\begin{eqnarray*}
\|C_i Q_g - Q_g \Lambda_i\|_F^2 &=& 
\|C_i  [Q_1 Q_2]   -  [Q_1 Q_2]
\begin{pmatrix}
  \Lambda_i\up{1} & 0 \cr  0 &   \Lambda_i\up{2} 
\end{pmatrix}  \|_F^2\\
&=& \| C_i Q_1 - Q_1 \Lambda_i\up{1} \|_F^2+\| C_i Q_2 - Q_2 \Lambda_i
\up{2}\|_F^2 . 
\end{eqnarray*} 
Clearly, the two terms on the right-hand side are independent.
What this means is that 
 we may minimize the sum (over $i$) of the first terms separately from the
sum of the second terms. The focus is on $Q_1$.
Alternatively,
we can  solve the same problem \nref{eq:obj} for 
diagonal $D_i$'s but now the matrix $Q$ is restricted to having
only $k$ columns. 

From what was just said, the matrix $Q_1$ obtained
from what we term a global optimum is the same as the matrix
$ Q $ obtained from minimizing 
\nref{eq:obj} with the additional constraint
of diagonality of the $D_i$'s. This is stated next.
\begin{proposition}
The matrix $Q_1$ extracted from the first $k$ columns of
the global solution $Q_g$ that minimizes 
\nref{eq:obj03}  is 
an optimal solution of Problem \nref{eq:obj} 
over orthogonal  $ n \times k$  matrices $Q$, and diagonal 
 matrices $D_i$, $i=1,\cdots, p$.
\end{proposition}

The proof of this result is straightforward and it is omitted.
Notice that we did not state the other half of the result, i.e., we did not
say  that   an optimal solution  $Q$
of \nref{eq:obj} with diagonal $D_i$'s, is essentially equal to the $Q_1$ 
extracted as shown above. This is due to non-uniqueness issues.
We cannot guarantee that the
solution $Q$  obtained by solving the reduced size problem \nref{eq:obj}
corresponds to the dominant subspace. For example, 
it can be that the algorithms being used  yield a lower value 
for the objective function, that is associated, say, with the
smallest eigenvalues. If it were possible to  guarantee that the
subspace $\Span{Q}$ corresponds to an approximate \emph{dominant}
 subspace  for all the
$C_i$'s, then we could  state that the solution are  identical with
that extracted from the global solution.
 However, this is harder to formulate rigorously because the
eigenspaces are not exact but only approximate. 
In practice, we will need to ensure that the eigenspace selected in 
whatever procedure is used to compute $Q$, will select the 
dominant eigenspace. 

\section{Computing an optimal $Q$}
In this section we discuss two distinct algorithms for computing
a matrix $Q$ that minimizes the objective function $\gamma$. 
The first one is based on the well-known subspace iteration algorithm 
and the second is  a gradient-based procedure.

\subsection{A subspace iteration procedure} \label{sec:subs} 
The simplest procedure for finding an optimal orthogonal 
$Q$ is to rely on the subspace iteration 
algorithm~\cite{Parlett-book,Saad-book3}. In the case of a single
matrix the algorithm amounts to taking some random matrix 
$X\in \RR^{n\times k}$ and computing the subspace spanned by $Y=A^j X$,
for some power $j$. The span of the columns of $Y$, is then
taken as an approximation to the dominant subspace. Since we have
$p$ matrices, we will proceed by combining the subspace iteration
with the SVD with a goal of extracting a `common dominant subspace'.

If we perform one single step of the subspace iteration,
i.e., taking $j=1$ in the above description, for each matrix $C_i$
then we would end-up with $p$ different subspaces of dimension $k$ each.
Putting these together yields a subspace of dimension $k\times p$.
The SVD can now be used to extract a nearly common subspace  of dimension
$k$ and the process is repeated. The algorithm is described next. 

\begin{algorithm}
  \caption{Subspace iteration}\label{alg:subsJ}
\begin{algorithmic}[1]
\State  Start : select initial $Q$ such $Q^T Q = I$
 \While{Not converged} 
 \For{ $j=1,\cdots, p$}
\State      Compute $X_j  = C_j Q$ 
 \EndFor 
\State    Let $X = [X_1, \cdots, X_p] $
\State  Compute $X = Q \Sigma V^T$  the SVD of $X$ 
\State  Define $Q := Q(:,1:k) $ [Matlab notation used]  
 \EndWhile
\end{algorithmic}
\end{algorithm}

Note  that when $p=1$ the algorithm amounts to a simple subspace iteration
algorithm with one matrix. Then, the SVD step in Lines~7-8 
replaces the Rayleigh Ritz process used in subspace iteration to 
extract the set of dominant eigenvectors.
In practice, the SVD calculations in Lines~7-8, can be performed using
Krylov subspace-type methods since  only the $k$ top 
left singular vectors of $X$ are sought.
One way to analyze then convergence of this algorithm is to view it from the 
perspective of a perturbed subspace iteration method~\cite{ys_subs14}.

\subsection{Gradient-based  Algorithms}\label{sec:grad}

It is natural to think of a gradient-type algorithm to maximize  the objective
function  \nref{eq:Psi} as  a means of balacing the objectives of invariance
and subspace targetting. 
This leads us to considering the Grassmannian alternative described at the end
of Section~\ref{sec:grass}. With a Grassmannian perspective the
gradient of $\psi$ is given by \nref{eq:gradPsi}.
With the notation $R_i \equiv R_i(Y)= C_i Q - Q D_{Y,i} $
this  gradient is written as: 
\eq{eq:grad1} 
G =   (1-\eta) \sum_{i=1}^p R_i + \eta \sum_{i=1}^p [ Y Y^T C_i R_i  + R_i D_{i,Y} ] 
\en
The next iterate will be of the form
\eq{eq:tilQ}
\tilde Q = Q + \mu G,
\en
where $\mu$ is to be determined.  As is known the 
direction of the gradient is a direction of increase for the objective
function $\psi$. It remains to determine how to select an
optimal $\mu$. First, we note that  $\mu $ cannot  arbitrarily large
because this would violate the restriction that the new $Q$ must be orthogonal.
In fact we must `correct' the non-orthogonality of the update \nref{eq:tilQ}.

Because $Q^T G = 0$ we have: 
\[
  \tilde Q^T \tilde Q = [Q + \mu G]^T [Q + \mu G]  = I + \mu^2 G^T G . 
\]
Let $G^T G = U \Sigma^2 U^T$ and define the diagonal matrix
\eq{eq:Psi3}
D_\mu \equiv [I + \mu^2 \Sigma^2]^{1/2} .
\en
In order to re-place
$\tilde Q$ on the Stiefel manifold without  changing its linear span
we will multiply it to the right by $U D_\mu \inv$, i.e., we define
\eq{eq:Qnew}
Q_{new} = \tilde Q U D_\mu\inv = (Q + \mu G) U D_\mu\inv . 
\en
With this we will have, 
\begin{eqnarray*}
  Q_{new}^T  Q_{new} 
 &=&  D_\mu\inv U^T [Q + \mu G]^T  \times [Q + \mu G] U  D_\mu\inv  \\
 &=&  D_\mu\inv U^T [I + \mu^2  G^T G ] U  D_\mu\inv  \\
  &=&  D_\mu\inv [I + \mu^2 \Sigma^2 ] D_\mu\inv    \\
  &=& I
\end{eqnarray*} 
as desired. If we set
\eq{eq:QuGu}
Q_u = Q U, \qquad G_u = GU ,
\en
the objective function $\psi$ for this $Q_{new}$  becomes the following function of
$\mu$:
\eq{eq:ratmu}
h(\mu) = \phi( [Q_u + \mu G_u] D_\mu\inv ) - \eta \gamma ([Q_u + \mu G_u] D_\mu\inv ) .
\en
This is a rational function which is asymptotic to a constant at infinity
and is an increasing function around $\mu = 0 $.
It is easy to devise procedures to \emph{approximately} maximize $h(\mu)$ by
sampling, e.g., uniformly.

\begin{algorithm}
  \caption{Gradient Ascent algorithm}\label{alg:gradA} 
  \begin{algorithmic}[1]
    \State  \textbf{Start:} Select initial $Q$ such that $Q^T Q = I$.

    \State     Compute $G$ from \nref{eq:grad1}.

    \While{$\| G \|_F > tol$}

\State  Compute and Diagonalize $G^T G$ as $G^T G = U \Sigma^2 U^T$

\State Compute $Q_u, G_u$ from \nref{eq:QuGu}.

\State Compute $\mu $ to approximately maximize \nref{eq:ratmu}

\State Set $ Q \ := \ (Q + \mu G) U [I + \mu^2 \Sigma^2]^{-1/2} $ .

\State     Compute $G$ from \nref{eq:grad1}. 
\EndWhile
\end{algorithmic}
\end{algorithm}

A gradient procedure may be appealing if a good approximate solution is already
known, in which case, the gradient algorithm may provide a less expensive
alternative to one step of the subspace iteration \ref{alg:subsJ}.  The
numerical experiments to be presented next will emphasize the subspace iteration
approach.

\section{Simulation results}
Since there is no joint  diagonalization algorithm adap\-ted  to high
dimensional data sets we will first  compare the proposed method  on 
small data set examples with a popular existing method and analyze its
behavior. In  this small  example, the
proposed method is compared to JADE \cite{Cardoso2}.
Second, the application of the proposed algorithm is extended
to a  high dimensional data  example. 

\subsection{Small diagonalizable data sets}
A  set of  $L$  diagonalizable  correlation matrices  of  size $N$  is
constructed using the  relation $C_{oD}^{l}=AC_{D}^{l}A^\top$. In this
case the  $N$ diagonal entries of  $C_{D}^{l}$ and the entries  of $A$
are independent  and identically distributed (i.i.d)  sampled from the
standard  normal distribution  and $C_{D}^{0}$  is taken to be the 
identity 
matrix. These  matrices can  be jointly  diagonalized with  any matrix
that differs from $A$ only by a permutation of the rows of $A^{-1}$ or
by scale factors multiplied to these rows. The matrix sizes were 
all multiples of 10 from $N=10$ to $N=100$, and we took $L=10$.

\subsection{Small approximately diagonalizable data sets}
The data set was generated the same way as above but each $C_{oD}^{l}$
is computed with an individual  matrix $A^{l}$. The entries of $A^{l}$
are obtained  as $a_{ij}^{l}=a_{ij}+n_{ij}$  where $n_{ij}$  are i.i.d
sampled from  $N(0,0.01)$. The matrix sizes $N$ were as above
($N=10,20,\cdots, 100$) and $L=10$.
For both data sets
the   methods   were   analyzed    using   the   performance   measure
\cite{Vollgraf2}
\begin{displaymath}
E=\frac{1}{N^{2}-N}\sum_{l=1}^{L}\parallel\textrm{Off}\left(Q^\top C_{D}^{l}Q\right)\parallel_{F}^{2}.
\end{displaymath}
In all cases the proposed method  and JADE were applied to 100 
generated data  sets and the  measure was  taken as the  mean obtained
over the  100 cases.  For the proposed approach  the dimension
$k$     was  fixed at   $k=5$     for     all    the     sizes
$N=10,20,\cdots,100$. Figure  1 illustrates  the results
for  the  diagonalizable data  sets.  The  figures correspond  to  the
boxplot of the  measure $E$ obtained over $100$ data  sets of the same
size but with a new realization of the matrix $A$. We can observe from
these  figures that  while $E$  increases with  the dimension  for the
method \cite{Cardoso2}, the  measure $E$ is relatively  stable for the
proposed  approach  as  $k=5$  for all  data  sets.  Furthermore,  the
proposed  method  outperform  \cite{Cardoso2}   in  term  of  $E$  and
computational time  as it took  0.066 sec  for the proposed  method to
joint diagonalize $10$ matrices of  size $N=100$ whereas it took 15.72
seconds for  JADE to  joint diagonalize  the same  data set.  Figure 2
illustrates  the variation  of the  $E$  as indicated  by the  boxplot
obtained over 100 diagonalizable  data sets of the  same size
$N=10$  but  with different  dimension  $k=3,4,5,6$  and $7$.  We  can
observe that once  the appropriate dimension is  selected the approach
is relatively stable as well. Figure  3 which is obtained based on the
small approximately diagonalizable data sets under the same conditions
(100 realization  of  the  measure   $E$)  illustrates  the
robustness of the proposed method to small errors.

\begin{figure}[!tbp]
  \centering
  \includegraphics[width=0.45\textwidth]{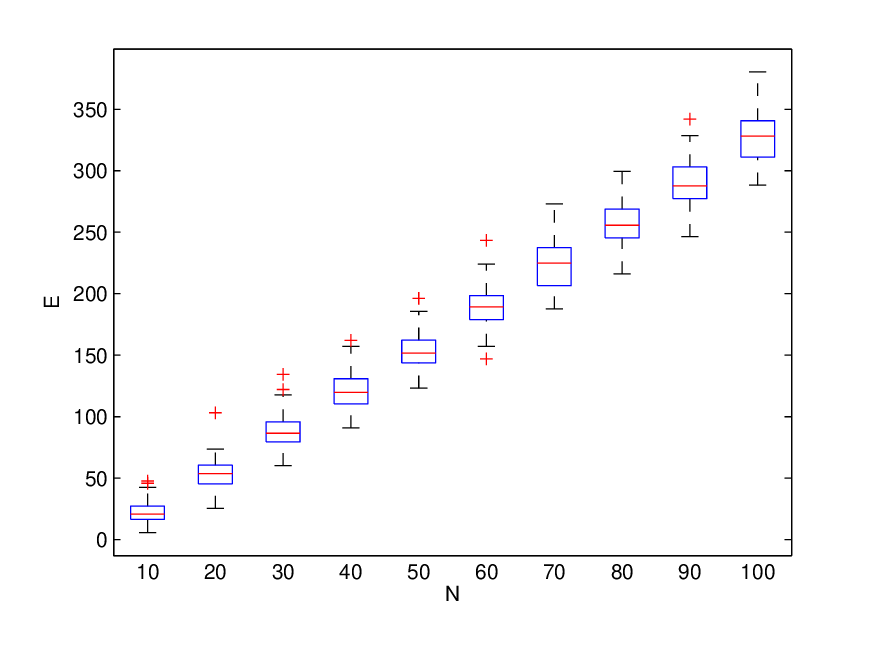}\label{fig:f1}
  \hfill
  \includegraphics[width=0.45\textwidth]{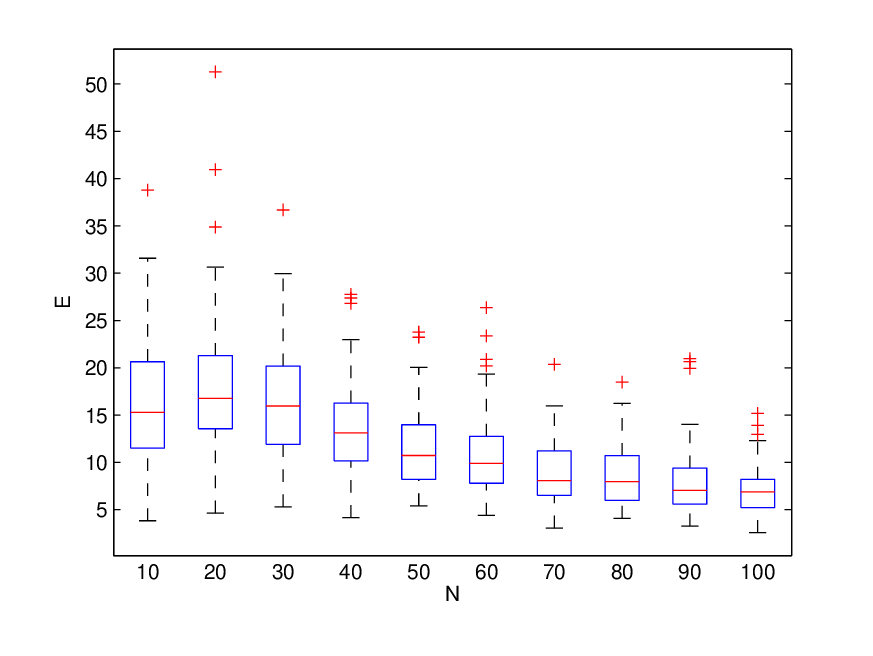}\label{fig:f2}
  \caption{Simulation results for the small diagonalizable data sets. The figures present the boxplots of the measure $E$ generated over 100 realization for \cite{Cardoso2} in the left figure and the proposed method in the right figure.}
\end{figure}

\begin{figure}[!tbp]
  \centering
  \includegraphics[width=0.45\textwidth]{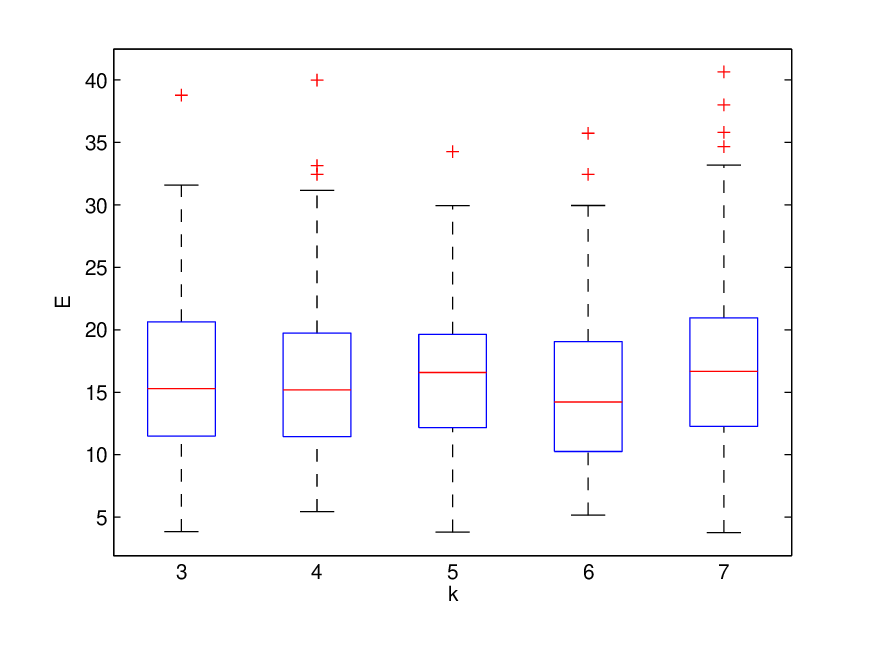}\label{fig:f3}
  \caption{Simulation results for the small diagonalizable data set of size $N=10$. The figures present the boxplots of the measure $E$ generated over 100 realization with the proposed method with dimension $k=3,4,5,6$ and $7$.}
\end{figure}

\begin{figure}[!tbp]
  \centering
  \includegraphics[width=0.45\textwidth]{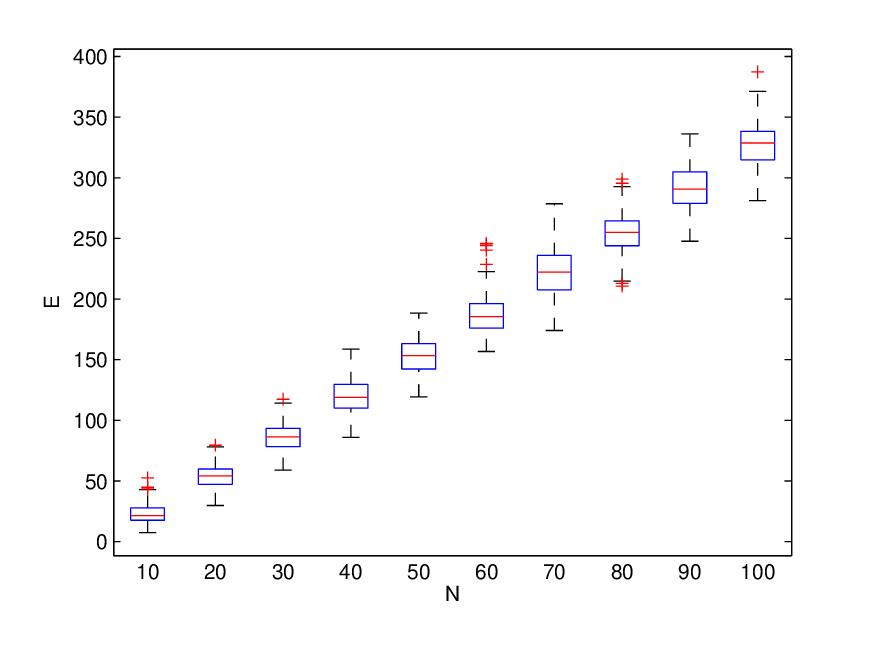}\label{fig:f4}
  \hfill
  \includegraphics[width=0.45\textwidth]{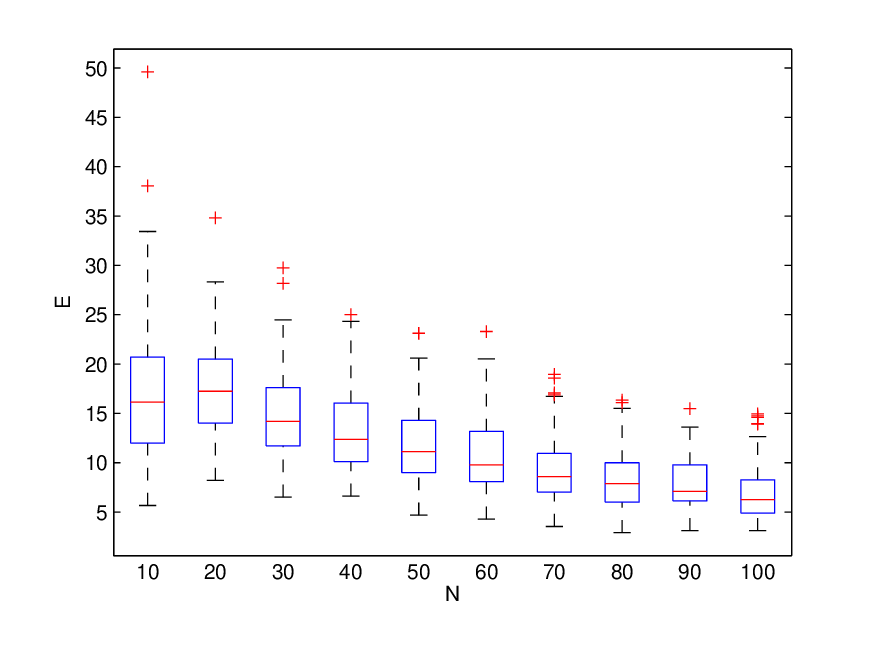}\label{fig:f5}
  \caption{Simulation results for the small approximately diagonalizable data sets. The figures present the boxplots of the measure $E$ generated over 100 realization for \cite{Cardoso2} in the left figure and the proposed method in the right figure.}
\end{figure}

\subsection{Restin state fMRI experiment}
In  this section  we evaluate  the performance  of the  proposed joint
matrices diagonalization algorithm on  a resting state fMRI experiment
data  set   \cite{Shehzad}.  Data-driven  methods   were  successfully
suggested and  applied to fMRI  data analysis. These  methods consider
the fMRI  time series measured at  each voxel as a  mixture of signals
localized  in   a  small  set   of  regions  and   other  simultaneous
time-varying  effects.  They isolate  the  spatial  brain activity  by
estimating a mixing  matrix and the sources that  define the spatially
localized neural dynamics.  Most data driven  fMRI analysis methods
use a  data matrix $\textbf{Y}$  formed by vectorizing   each time
series  observed in  every  voxel creating  a  matrix $\textbf{Y}$  of
dimension $N\times n$  where $n$ is the number of  time points and $N$
the number of voxels, $N\gg n$.\\
These methods consider $\textbf{Y}$ as
the  mixture  and  factorize  it   into  latent  sources  through  the
decomposition   into  matrices   $\textbf{Y}=\textbf{A}\textbf{X}$,  a
mixing     matrix     $\textbf{A}$     and     a     source     matrix
$\textbf{X}$.  Data-driven methods  are suitable  for the  analysis of
fMRI data as they minimize the assumptions on the underlying structure
of the  problem by  decomposing the  observed data  based on  a factor
model and  a specific  constraint. Different  constraints have  led to
different  data-driven  methods.  For example,  the  maximum  variance
constraint   has   led   to   principal   component   analysis   (PCA)
\cite{Firston},  the independence  constraint has  led to  independent
component analysis  (ICA) \cite{McKeown}  and sparsity  constraint has
led to dictionary  learning \cite{Lee}.

Recently, ICA  has become a
widespread  data-driven  method  for  fMRI analysis.  It  has  led  to
temporal ICA  (tICA, for the format  of the data described  above) and
spatial ICA (sICA)  \cite{McKeown}. In this experiment  we applied the
proposed joint diagonalization approach on twenty correlation matrices
of size  $939\times 939$ obtained from  a data set of  size $939\times
197$. This data  set was constructed from the slice  41, which we know
contains  the activated  regions  of the  default  mode network  (DMN)
\cite{Shehzad}.  For comparison  we  used the  tICA  approach. We  can
observe  from  figure 4  that  both  the  proposed approach  and  tICA
recovered the  connected regions of  the DMN; the  posterior cingulate
cortex  (PCC), medial  pre-frontal  cortex (MFC),  and right  inferior
parietal lobe  (IPL). Since there is  no gold standard reference  for DMN
connectivity  available,   we  relied on  the similarity  of
temporal   dynamics  of   DMN  based   modulation  profile   with  PCC
representative   time-series.   The   similarity  measure   used   was
correlation and it was estimated as $>0.75$ for all the algorithms.

\begin{figure}[!tbp]
  \centering
  \includegraphics[width=0.45\textwidth]{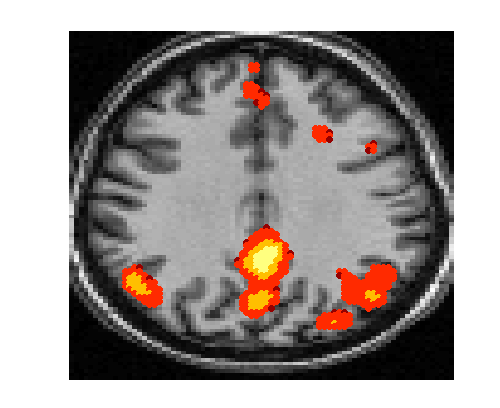}\label{fig:f6}
  \hfill
  \includegraphics[width=0.45\textwidth]{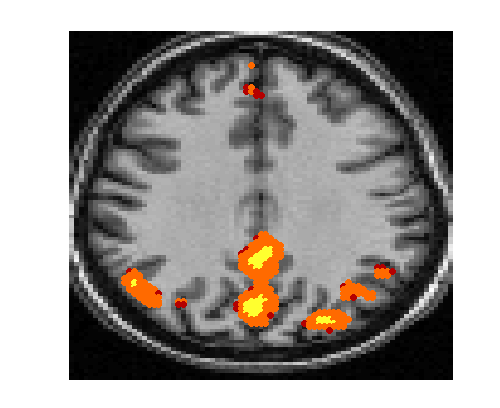}\label{fig:f7}
  \caption{Resting state fMRI results. Left: tICA \cite{McKeown}, right: proposed approach}
\end{figure}

\section{Conclusion}
This paper discussed what might be termed a natural extension of the
problem of joint diagonalization for the situation when the matrices
under consideration are too large for standard algorithms to be applied. 
This extension consists of solving the problem in a subspace of small
dimension, leading  to the minimization of an objective function to
produce  an orthonormal basis of the desired subspace.
A few theoretical results have been shown that characterize this
optimum and establish a few equivalences. An algorithm based
on a variant of subspace iteration was proposed to solve the
problem and was tested on a few examples.
One issue that still remains to be addressed is to show that the joint
approximate eigenspace 
to which the algorithm converges corresponds to a dominant 
eigenspace for all the $C_i$'s.
This requirement is somewhat difficult to 
formulate rigorously due to the approximate nature of the eigenspaces 
involved. In the easiest case of exact joint 
joint diagonalization, we would have $C_i Q - Q D_i =0$ for all $i$ and
we would ask that in addition $Q$ be an eigenspace associated with the $k$ 
dominant eigenvalues for each $C_i$. In the approximate case,
$\| C_i Q - Q D_i \| $ can only be required to be small.
Then, because $Q$ is only an approximate eigenspace, it is meaningless
to  demand that the associated approximate eigenvalues
be the largest ones for each $C_i$.
Although the subspace iteration algorithm is likely to deliver 
a solution that more or less satisfies this requirement, 
the theoretical foundation, as well as a rigorous formulation of the 
result, are yet to be established. 

\bibliographystyle{siam}
\bibliography{local,bibliography}

\end{document}